\documentclass{amsart}
\usepackage{amsmath, amssymb}
\usepackage{fullpage}
\usepackage{graphicx}

\newenvironment{example}{
\begin{quotation}\noindent{\bfseries Example.\/} }{\end{quotation}}

\def\EE{\ensuremath{ E}}

\newcommand{\rsc}[2]{\left\langle#1, #2\right\rangle} 
\def\field{K}

\newcommand\exq[1]{\ifthenelse{\equal{#1}{}}{\operatorname{ex}_q}{\operatorname{
ex}_q\left(#1\right)}}

\newcommand\bN{\mathbb N}

\newcommand\bZ{\mathbb Z}

\def\sF{\ensuremath{\mathsf{F}} } 
\def\sP{\ensuremath{\mathsf{P}} } 
\def\sE{\ensuremath{\mathsf{E}} }
\def\sG{\ensuremath{\mathsf{G}} }

\def\sX{\ensuremath{\mathsf{X}} }
\def\sL{\ensuremath{\mathsf{L}} }
\def\sC{\ensuremath{\mathsf{C}} }

\def\C{\ensuremath{\mathbb C}}

\def\N{\ensuremath{\mathbb N}}

\def\spa{{\sf scalar\_de}}
\def\ham{{\sf hammond}}

\usepackage{amssymb,amsmath, amsthm}
\newtheorem{thm}{Theorem}[section]
\newtheorem{cor}[thm]{Corollary}

\newtheorem{defn}{Definition }
\newenvironment{pf}{\begin{proof}}{\end{proof}}

\title{Automatic enumeration of regular objects}
\author{Marni Mishna}
\email{mmishna@sfu.ca}
\address{Department of Mathematics\\ Simon Fraser University\\8888 University Drive\\ Burnaby,
  Canada\\V5A 1S6} 
\thanks{The author thanks the Canadian Natural Science and
  Engineering Research Council for funding this work through the PDF program.}
\thanks{This work was completed while at  LaBRI, Universit\'e 
  Bordeaux I, and at the Fields Institute, Toronto, Canada.}
\thanks{{\bf keywords:} Asymptotic enumeration, automatic combinatorics, generating
  functions, symmetric functions}
\thanks{{\bf Mathematics Subject Classifications:} 05A16, 05C30}

\begin{document}

\begin{abstract}
  We describe a framework for systematic enumeration of families
  combinatorial structures which possess a certain regularity. More
  precisely, we describe how to obtain the differential equations
  satisfied by their generating series. These differential equations
  are then used to determine the initial terms in the counting sequence and for
  asymptotic analysis.  The key tool is the scalar product for
  symmetric functions.\\

\end{abstract}
\maketitle
\section*{Introduction}
Some classes of combinatorial objects naturally possess a substantial
amount of symmetry and when formal sums of monomials encoding
some parameter of interest are taken over the entire class, symmetric
functions, or symmetric series appear. There has been
some recent activity to determine how to extract enumerative series of
sparse sub-families of these classes directly from the symmetric
functions. The principle can be illustrated with one well studied
example, the subset of labelled graphs in which the degree of each
vertex is a fixed value, say $k$, known as the {\em $k$-regular
  graphs}. Here, we encode a graph by its degree sequence. When
we consider the sum of this encoding over all graphs, they are encoded
by the infinite product
\begin{equation}\label{eq:allgraphs}
  G(x_1, x_2, \ldots)=\prod_{i<j}(1+x_ix_j).
\end{equation}
This is a well known symmetric series. Further, we remark that the
coefficient of $x_1^kx_2^k\cdots x_n^k$ in the formal power series
development gives the number of  labelled $k$-regular
graphs on $n$ vertices.

Such a coefficient extraction can be set up
as a multidimensional Cauchy integral, as described by McKay for regular
tournaments and Eulerian digraphs~\cite{McKay90}, and by McKay and Wormald
for graphs with a fixed degree sequence~\cite{McWo90}. However, in
general this may not be a useful or practical formula. 

Indeed these techniques are well developed, and provide general
formulae which we cannot currently obtain with the methods here, but they are difficult to make systematic,
as they contain a saddle point analysis to make the asymptotic
estimate which may be quite fine and specific to the problem.

The primary goal here is to lucidly illustrate how techniques for computing the scalar product
of symmetric functions in~\cite{ChMiSa05} can be a part of an
essentially algorithmic process for asymptotic analysis. At the heart of the method is
the fact that the scalar product of symmetric functions preserves a
notion of D-finiteness~\cite{Gessel90}, and, thanks to the algorithms in~\cite{ChMiSa05}, this
result is effective. 

We begin with a short recollection of symmetric series and
D-finiteness, and a brief discussion on some places that D-finite
symmetric series appear in combinatorics.  We analyse graphs with
fixed finite degree sets, and hypergraphs. Finally, in
Section~\ref{sec:asympt}, we have the results of our semi-automated
asymptotic analysis of these classes.

\section{Symmetric series and D-finiteness}
We provide a basic summary of symmetric functions in order to
establish notation. The reader is directed
to MacDonald's book~\cite{Macdonald95} for full details.

Denote by $\lambda=(\lambda_1,\dots,\lambda_k)$ a partition of the
integer~$n$.  This means that $n=\lambda_1+\dots+\lambda_k$
and~$\lambda_1\ge\dots\ge\lambda_k>0$, which we also
denote~$\lambda\vdash n$.  Partitions serve as indices for the five
principal symmetric function families that we use: homogeneous
($h_\lambda$), power ($p_\lambda$), monomial ($m_\lambda$), elementary
($e_\lambda$), and Schur ($s_\lambda$).  These are series in the
infinite set of variables, $x_1,x_2,\dotsc$ over a field~$\field$ of
characteristic~0.  When the indices are restricted to all partitions
of the same positive integer~$n$, any of the five families forms a
basis for the vector space of symmetric polynomials of degree~$n$
in~$x_1,x_2,\dotsc$\@{}  On the other hand, the family of $p_i$'s indexed
by the integers~$i\in\bN$ generates the algebra~$\Lambda$ of symmetric
functions over~$\field$: $\Lambda=\field[p_1, p_2,\dotsc]$.
Furthermore, the~$p_i$'s are algebraically independent over~$\bZ$.

Generating series of symmetric functions live in the larger ring of
symmetric series, $\field[t][[p_1,p_2,\dotsc]]$. There, we have the
generating series of homogeneous and elementary functions:
\[
H(t)=\sum_n h_n t^n =\exp\left(\sum_i p_i\frac{t^i}i\right),\qquad
E(t)=\sum_n e_n t^n =\exp\left(\sum_i (-1)^ip_i\frac{t^i}i\right).
\]
We often refer to $H=H(1)$ and $E=E(1)$.

Alternatively, the power notation $\lambda=1^{n_1}\dotsm k^{n_k}$ for
partitions indicates that $i$~occurs $n_i$~times in~$\lambda$,
for~$i=1,2,\dots,k$.  The normalization constant
\[
z_\lambda:= 1^{n_1}n_1!\dotsm k^{n_k}n_k!
\]
plays the role of the square of a norm of~$p_\lambda$ in the following important
formula:
\begin{equation}\label{eq:scalp}
\rsc{p_\lambda}{p_\mu}=\delta_{\lambda,\mu}z_\lambda,
\end{equation}
where $\delta_{\lambda,\mu}$~is~1 if $\lambda=\mu$ and 0 otherwise. 

The scalar product is a basic tool for coefficient extraction. Indeed,
if we write $F(x_1,x_2,\dotsc)$ in the form~$\sum_\lambda f_\lambda
m_\lambda$, then the coefficient of $x_1^{\lambda_1}\dotsm
x_k^{\lambda_k}$ in~$F$ is $f_\lambda=\rsc F{h_\lambda}$.  Moreover,
when $\lambda=1^n$, the identity $h_{1^n}=p_{1^n}$ yields a simple way
to compute this coefficient when $F$ is written in the basis of
the~$p$'s. When viewed at the level of generating series, this fact
gives the following theorem:
\begin{thm}[Gessel\cite{Gessel90}; Goulden \& Jackson\cite{GoJaRe83}]\label{thm:theta}
  Let $\theta$ be the $\field$-algebra homomorphism from the algebra
  of symmetric functions over~$\field$ to the algebra~$\field[[t]]$ of
  formal power series in~$t$ defined by $\theta(p_1)=t$,
  $\theta(p_n)=0$ for $n>1$. Then if $F$ is a symmetric function,
\[\theta(F)=\sum_{n=0}^\infty a_n\frac{t^n}{n!},\]
where $a_n$ is the coefficient of $x_1\dotsm x_n$ in~$F$.
\end{thm}

To end our brief recollections of symmetric functions recall that
plethysm is a way to compose symmetric functions.  An inner law
of~$\Lambda$, denoted $u[v]$ for $u,v$ in~$\Lambda$, it satisfies the
following rules~\cite{Stanley99}, with $u, v, w\in\Lambda$ and~
$\alpha,\beta$ in~$\field$
\begin{equation*}
  (\alpha u+\beta v)[w]=\alpha u[w]+\beta v[w],
  \quad
  (uv)[w]=u[w]v[w],
\end{equation*}
and if  $w=\sum_\lambda c_\lambda p_\lambda$ then
$p_n[w]=\sum_\lambda c_\lambda p_{(n\lambda_1)}p_{(n\lambda_2)}\ldots$.
For example, consider
that $w[p_n]=p_n[w]$, and in particular that $p_n[p_m]=p_{nm}$.  In a
mnemonic way:
\[
w[p_n]=w(p_{1n},p_{2n},\dots,p_{kn},\ldots)
\qquad\text{whenever}\qquad
w=w(p_1,p_2,\dots,p_k,\ldots).
\]

\subsection{D-finite multivariate series}
Recall that a series $F\in\field[[x_1,\dots,x_n]]$ is {\em
D-finite\/} in $x_1,\dots,x_n$ when the set of all partial derivatives
and their iterates, $\partial^{i_1+\dots+i_n}F/\partial x_1^{i_1}\dotsm
\partial x_n^{i_n}$, spans a finite-dimensional vector space over the
field $\field(x_1,\dots, x_n)$.  A {\em D-finite description\/} of a
series~$F$ is a set of differential equations which
establishes this property.  A typical example of such a set is a
system of $n$~differential equations of the form
\[
q_1(x)f(x)+q_2(x)\frac{\partial f}{\partial x_i}(x)+\dots+
q_k(x)\frac{\partial^kf}{\partial x_i^k}(x)=0,
\]
where $i$~ranges over $1,\dots,n$, each~$q_j$ is in $\field(x_1,\dots,x_n)$
for $1\leq j\leq k$, and $k$ and~$q_j$ depend on~$i$.

Such a system is a typical example of a D-finite description of a
functions, and often this will be the preferred form for
manipulating~$f$. In truth we can accept any basis which generates the
vector space of partial derivatives, but in the applications below,
this form is particularly easy to obtain.

\subsection{D-finite symmetric series}
The following definition of D-finiteness of series in an infinite
number of variables is given by Gessel~\cite{Gessel90}, who had 
symmetric functions in mind. A series $F\in\field[[x_1,x_2,\dotsc]]$
is {\em D-finite\/} in the $x_i$ if the specialization to 0 of all but
a finite (arbritrary) choice, $S$, of the variable set results in a D-finite
function (in the finite sense). In this case,
many of the properties of the finite multivariate case hold true. One
exception is closure under algebraic substitution, which requires
additional hypotheses.

The definition is then tailored to symmetric series by considering the
algebra of symmetric series as generated over~$\field$ by the set
$\{p_1,p_2,\dotsc\}$: a symmetric series is called {\em D-finite\/}
when it is D-finite in the $p_i$'s\footnote{This is interestingly
  enough {\em not\/} equivalent to D-finiteness with respect to either the $h$ or
  $e$ basis.}.
\begin{example}
Both $H(t)$ and $E(t)$ are D-finite symmetric functions, as for any
specialization of all but a finite number of the $p_i$'s to 0 results
in an exponential of a polynomial. Similarly, $\exp(h_kt)$ is D-finite
because $h_k=\sum_{\lambda\vdash k} p_\lambda$ is a polynomial in the $p_i$s. 
\end{example}

The closure under Hadamard product of D-finite
series~\cite{Lipshitz88} yields the consequence:
\begin{thm}[Gessel]\label{thm:rsc_pres_df}
Let $f$ and $g$ be elements of
$\field[t_1,\dots,t_k][[p_1,p_2,\dotsc]]$, D-finite in the $p_i$'s
and $t_j$'s, and suppose that $g$ involves only finitely many of the
$p_i$'s. Then $\rsc fg$ is D-finite in the $t_j$'s provided it is
well-defined as a power series.
\end{thm}

\subsection{Effective calculation and algorithms}
In our initial study~\cite{ChMiSa05} we gave an algorithm
which, given a D-finite descriptions of two functions satisfying the
hypothesis of Theorem~\ref{thm:rsc_pres_df}, determines a D-finite
description of the series of the scalar product. Henceforth, we shall
refer to this algorithm as~\spa. As we noted in~\cite{ChMiSa05}, a
second algorithm,~\ham, based on the work of Goulden, Jackson and
Reilly~\cite{GoJaRe83} applies in the case when $g=
\exp(h_nt)$, which we shall see is precisely how one can
extract the exponential generating series of sub-classes with
``regularity''.  They are implemented in Maple, are are
available for public distribution at the website {\tt
  http://www.math.sfu.ca/\verb+~+mmishna}. Maple worksheets illustrating the
calculations presented are also available at that same site. 

\section{D-finite symmetric series appear naturally in combinatorics}
Species theory (in the sense of~\cite{BeLaLe98,Joyal81}) is formalism
for defining and manipulating combinatorial structures which relates
classes to encoding series. An important connection to our work here
is that the series for structures we consider are D-finite symmetric
series, and many of the natural combinatorial actions preserve
D-finiteness on the level of these series.

The reader unfamiliar with species is heartily encouraged to
consult~\cite{BeLaLe98}. A {\em species\/} associates to every set a
family of structures in a way such that two sets of the same
cardinality yield the same family, upto isomorphism. For example, the
species of sets $\sE$ on the underlying set $U$ is simply
$\sE[U]=U$. The species of lists $\sL[U]=\{(x_1, x_2, \dots, x_n):
  x_i\in U, n=0,1,2\dots\}$, is the set of finite ordered collections
  of elements.  The atomic species, $X[U]$ is $U$ if $U$ contains a
  single element, and is empty otherwise.

The theory of species develops a rigorous
formalism which allows a sort of calculus of combinatorial
families. For example we construct lists of length 4 from our atomic
species via multiplication: $\sL_4[U]=X^4[U]$.

The key feature that we use are that for every combinatorial family
(species) $\sF$ that one can define, there is an associated cycle
index series $Z_\sF$ and an asymmetric cycle index series~$\Gamma_\sF$
both of which are symmetric series. Recall for any species \sF its
cycle index series $Z_\sF$ is the series in $\C[\![p_1, p_2,
\ldots]\!]$ given by
\begin{equation}\label{eq:cycle_index}
Z_\sF(p_1, p_2, \ldots)
 :=\sum_n\sum\limits_{\lambda\vdash n}\operatorname{Fix} \sF[\lambda]
   \frac{p_1^{m_1}p_2^{m_2}\cdots p_k^{m_k}}{z_\lambda},
   \end{equation}
   where the value of $\operatorname{Fix}\sF[\lambda]$ is the
   number of structures of $\sF$ which remain fixed under some labelling
   permutation of type\footnote{A
     permutation of type $(1^{m_1}, 2^{m_2}, \ldots)$ has $m_1$ fixed
     points, $m_2$ cycles of length 2, etc.} $\lambda$, and $m_k$ gives
   the number of parts of $\lambda$ equal to~$k$. 

   The definition of the {\em asymmetry index series} of a
   species~$\sF$, denoted~$\Gamma_\sF$, as introduced by
   Labelle~\cite{BeLaLe98} is related, but more subtle. The
   series~$\Gamma$ behaves analytically in much the same way as the
   cycle index series, notably, substitution (in almost all cases) is
   reflected by plethysm, etc. Essentially, this series counts the objects with no
   internal symmetry.  Table~\ref{tab:smallspec} contains some small
   examples of both series.

\begin{table}[t]
\center
\begin{tabular}{lll|lll}\small
Object&Series&Symmetric function&Object&Series&Symmetric function\\\hline\hline
2-sets&$\Gamma_{\sE_2}$&$e_2=\frac{p_1^2}{2}-\frac{p_2}{2}$&2-multisets&$Z_{\sE_2}$&$h_2=\frac{p_1^2}{2}
+\frac{p_2}{2}$\\
3-sets&$\Gamma_{\sE_3}$&$e_3$           &3-multisets&$Z_{\sE_3}$&$h_3$\\ 
4-sets&$\Gamma_{\sE_4}$&$e_4$           &4-multisets&$Z_{\sE_4}$&$h_4$\\
$k$-sets &$\Gamma_{\sE_k}$&$e_k$         &$k$-multisets&$Z_{\sE_k}$& $h_k $\\
3-cycles &$Z_{\sC_3}$&$\frac{p_1^3}{3}+\frac{p_3}{3}$ & triples &$Z_{X^3}$&$p_1^3 $\\
4-cycles &$Z_{\sC_4}$&$\frac{p_1^4}{4}+\frac{p_2^2}{12}+\frac{p_4}{12}$& 4-arrays&$Z_{X^4}$&$p_1^4 $\\
5-cycles &$Z_{\sC_5}$&$\frac{p_1^5}{5}+\frac{p_5}{30}$ & 5-arrays& $Z_{X^5}$&$p_1^5 $\\
$k$-cycles &$Z_{\sC_k}$& $\sum_{cd=k} \phi(d)\frac{p_d^c}{k!}$ &
$k$-arrays&$Z_{X^k}$& $p_i^k$\\\hline\hline\\
\end{tabular}\\

\caption{Index series of small species and their corresponding
  symmetric functions}
\label{tab:smallspec}
\hrule
\end{table}

In a fashion similar to the cycle index series,~$\Gamma_\sF$ arises
through the enumeration of colourings of asymmetric $\sF$-structures.

A notable example is the species of sets, $\sE$. Recall for any finite
set $U$ we have that $\sE[U]=U$. The two series above turn out to be
$Z_\sE=\exp(\sum_n p_n/n)=\sum_n h_n$
and~$\Gamma_\sE=\exp(\sum_n(-1)^n p_n/n)=\sum_n e_n$. 

The primary advantage of this approach, as is true with any generating
series approach, is that natural combinatorial operations (set,
cartesian product, substitution) coincide with straighforward analytic
operations (sum, product, plethystic substitution)\footnote{This is
 slightly less true with the asymmetry index series, but true enough for our purposes.}.

The {\em exponential generating series\/} of a species $\sF$ is the
sum $\sF(t)=\sum_n |\sF[n]| \frac{t^n}{n!}$, where $|\sF[n]|$ is the number of
structures of type $\sF$ on a set of size
$n$. The {\em ordinary generating
function}, $\widetilde{\sF}(t)$, is the sum $\sF(t)=\sum_n
\operatorname{Orb}(\sF[n]) t^n$, where $\operatorname{Orb}(\sF[n])$ is
the number structures of $\sF$ on a set of size $n$ distinct up to
relabelling. Also recall the notation $[x^n]f(x)$ refers to the
coefficient of $[x^n]$ in the expansion of $f(x)$. This definition
extends likewise to monomials.  

The next result is essentially a collection of known results and basic facts of
D-finite series. 
 \begin{thm}\label{thm:Dfin_species}
    Suppose $\sF$ is a  species such that $Z_\sF$ is a
    D-finite symmetric series and write $p_n=x_1^n+x_2^n+\ldots$. Then
    all of the following series are D-finite with respect to $t$:
  \begin{enumerate}
  \item The exponential generating function $\sF(t)$;
  \item The ordinary generating
  function $\widetilde{\sF}(t)$, if the additional condition that $Z_\sF(p_1, p_2, \ldots)$ is D-finite with
  respect to the $x_i$ variables is also true;
  \item The series $\sum_n \left([x_1^k\cdots
      x_n^k]Z_\sF\right) \frac{t^n}{n!}$, for fixed $k$;
  \item  The series 
$\sum_n\sum_{\bar k\in S^n}\left(\left[ x_1^{k_1}x_2^{k_2}\dots
    x_n^{k_n}\right] Z_\sF\right) t^n/n!$, for any finite set $S\subset \bN$.
  \end{enumerate}
\end{thm}

 \begin{pf}
The first two parts are proved using two basic results about cycle
index series:
\begin{equation*}
  \sF(t)=Z_\sF(t, 0, 0, \ldots) \quad \text{ and } \quad
  \widetilde{\sF}(t)=Z_\sF(t, t^2, t^3, \ldots)
\end{equation*}
The first specialization is well-known to preserve
D-finiteness~\cite{Stanley80} for any $n$. 
The additional condition on the second item is sufficient to  prove
the D-finiteness since the stated substitution is the same as
$x_1\mapsto t$, and $x_i\mapsto 0$ otherwise. 

The third item of the proposition is proved by the expression 
\begin{equation*}
  \sum_n \left([x_1^k\cdots x_n^k]Z_\sF\right) \frac{t^n}{n!}=\rsc{Z_\sF}{\exp(th_k)},
\end{equation*}
which is D-finite by Theorem~\ref{thm:rsc_pres_df}.

The final item  of the proposition is true because the series is equal to
\begin{equation*}
 \rsc{Z_\sF}{\exp(t\sum_{i\in S}h_i)},
\end{equation*}
which is also D-finite by Theorem~\ref{thm:rsc_pres_df}.
\end{pf}

We have one large class of species for which the cycle index series is
D-finite. All of our examples come from this class.
\begin{thm}
Let $\sE$ be the species of sets and let~$\sP$ be a polynomial
species with finite support\footnote{Species which
can be written as polynomials of {\em molecular
  species}. For example, every species in Table~\ref{tab:smallspec} is
polynomial.}. Then~$\sF=\sE\circ \sP$ describes a species for
which~$Z_\sF$ is a D-finite symmetric series and provided that $\sP(0)=0$,~$\Gamma_\sF$ is also a D-finite symmetric series.
\end{thm}
\begin{proof}
If~$\sP$ is a polynomial species, then its cycle index series is
a polynomial in the $p_i$'s, say $P(p_1, \ldots, p_n)$. Composition of
species is reflected in the cycle index series by plethysm, thus
\[
Z_\sF=\exp(\sum_k p_k/k)[P(p_1, \ldots, p_n)]=\exp\left(\sum P(p_k,
p_{2k},\dots, p_{nk})/k\right).
\]
For any specialization of all but a finite number of $p_i$ to zero,
this gives an exponential of a polynomial, which is clearly
D-finite. Thus, $Z_\sF$ is D-finite. We can similarly show that
$\Gamma_\sF$ is also D-finite under the stated conditions, since the
composition also results in a plethystic composition.
\end{proof}

Our concluding remarks in Section~\ref{sec:conclusion} address the
more general question of combinatorial criteria on a species $\sF$
that ensure that $Z_\sF$ or $\Gamma_\sF$ are D-finite.
\section{Using species to describe regular graph-like structures}
Ultimately our goal is to generalize the well-studied case of $k$-regular
graphs to other structures whose cycle index series are D-finite. To
do so, we express the graph encoding by degree sequence as symmetric
series, and describe how to find such a representation in general
using species theory.

In Eq.~\eqref{eq:allgraphs} we define $G(x_1, x_2, \ldots)$ as
the encoding over all graphs of their degree sequence and we express
this as an infinite product. It turns out that this series is
equivalent to $E[e_2]$, which is equivalent to $\Gamma_{\sE\circ
  \sE_2}$. The equivalent series for multigraphs (with loops) is equal
to $H[h_2]=Z_{\sE\circ\sE_2}$, and thus suspecting an explanation via
species, we investigate this connection. Specifically,  how do we
construct symmetric function equations to describe the generating
functions of different families of objects, such as hypergraphs.

We begin with the remark that $\sF=\sE\circ\sE_2$ is {\em not\/} the
species of graphs. It is the species of partitions into 2-sets. For
example, $\{\{1,4\},\{2,6\}, \{3,7\},\{5,8\}\}$ is an element
of~$\sF$, and we should not think that this is the graph on 8
vertices, with four edges, rather it just gives the basic structure,
i.e. four edges.

We express the Polya cycle index in the power series
symmetric function. As a series in the symmetric $x_i$
indeterminates, it is an inventory of distinct (non-isomorphic)
colourings of the elements of the species. 
  For example, the non-isomorphic colourings (by positive
integers, say) of the set  $\{a,b\}$ is the set of maps
$\{(a,b)\mapsto(i,j)\in\bN^2: i\leq j\}$, and the inventory of all
such colourings is $\sum_{i\leq j} x_ix_j=h_2$.

\begin{figure}
\center
\includegraphics[width=5cm]{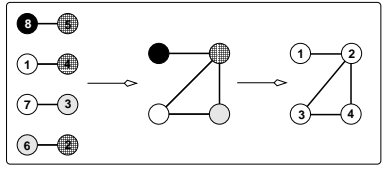}
\caption{The graph associated to the coloured set partition $\{ \{1_3,4_2\}, \{2_2,6_4\},\{3_4,7_3\},\{5_2,8_1\}\}$}
\label{fig:graph}
\hrule
\end{figure}

A colouring of an element in $\sE\circ\sE_2$ gives rise to a graph
(See Figure~\ref{fig:graph}) and two colourings are isomorphic if one
is a graph relabelling of the other. The monomial encoding a colouring indicates
how many time each colour was used, that is, in how many edges the
colour appears, that is, the degree of the vertex represented by the
colour.
 
We restate this correspondence. The species $\sE\circ\sE_2$ indicates
the structure-- sets of pairs. The cycle index series
$Z_{\sE\circ\sE_2}$ encodes non-isomorphic colourings of elements,
which are in turn equivalent to labelled multi-graphs. The sets of pairs
indicate edges, and the colours indicate vertices.

For many applications, like regular graphs, we would like to count
colourings without repetition. In this case, we do not allow
repetition of a colour in a given object, hence to encode a $k$-set,
each colour appears exactly once, and this is precisely the notion of
asymmetry in the asymmetry cycle index, and thus we use $\Gamma$
instead of $Z$. Remark,
and thus
\begin{equation*}
  \Gamma_{\sE_k} = e_k \quad\text{and thus, } \quad\Gamma_\sE=E.
\end{equation*}
Taking the same species $\sE\circ\sE_2$ as above, and using the
asymmetry index series with a similar argument, we get that
$\Gamma_{\sE\circ\sE_2}=\EE[e_2]$ encodes simple graphs without loops
on the set of colours precisely as is determined by
Eq.~\eqref{eq:allgraphs}: $E[e_2]=\prod_{i<j}(1+x_ix_j)$.  This gives
us a way to have {\em direct access to monomial encodings of
  combinatorial objects\/}, as symmetric functions expressed in common
bases, like the power sum basis.  These two series are compatible and,
one can show that graphs with loops are encoded by~$E[h_2]$, and
graphs with multiple edges, but no loops are given by~$H[e_2]$.

\begin{figure}
\center
\includegraphics[height=3.5cm]{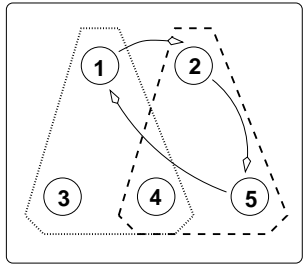}
\caption{A structure composed of a set of smaller structures (cycles
  and 3-sets).}
\label{fig:multi}
\hrule
\end{figure}

More generally, we can consider any structure which is built as a set
of objects from a finite set of classes. Figure~\ref{fig:multi} shows
a more general object built as a set of cycles and sets.  In this
framework it is encoded by the monomial~$x_1^2x_2^2x_3x_4^2x_5^2$, and
thus we see that regularity in this situation refers to the number of
times each label appears in one of the smaller substructures.

Using this framework we can examine other species of structures built
up from smaller objects. These species are such that both $Z_\sF$
and~$\Gamma_\sF$ give rise to interesting combinatorial objects.

We can produce enumerative results for objects all of the same
flavour: labelled sets of objects in which
there is a certain regularity. We begin with a natural generalization
of $k$-regular graphs, and then we consider other types of
objects such as hyper-graphs, and directed graphs.

\subsection{$S$-regular graphs}
A graph is {\em $S$-regular} if the set of vertex degrees in the graph
is a subset of $S$. For example, a graph is $\{i,j\}$-regular if every vertex
is of degree $i$ or $j$.

It does not seem that the asymptotic enumeration of these objects has been
directly considered before. It is, in some sense, a variation of the 
asymptotic number of labelled graphs with a given degree sequence,
which has been considered by Bender and Canfield~\cite{BeCa78} and
McKay and Wormald~\cite{McWo90}, and may very well be computable from
this. 


Thus, the scalar product which represents the generating
series for the number of $\{i, j\}$-regular graphs is given by
incorporating this factor, which ultimately greatly simplifies the
calculation. The exponential generating series for the number of
$\{i,j\}$-regular graphs is given by
\begin{align*}
G_{i,j}(t)&=
\big\langle E[e_2], \exp\left(t(h_i+h_j)\right) 
  \big\rangle.
\end{align*}

This is clearly D-finite, and computable using \spa~(although
not~\ham). Furthermore, by a similar computation, we have the
following result.
\begin{thm} 
  The number of~$S$-regular graphs is D-finite for any
  finite~$S\subset\N$, and its exponential generating series is given
  by the scalar product of symmetric functions,
\begin{equation*}
  G_S(t)=\left\langle E[e_2], \exp\left(\sum_{i\in S} h_it\right)
  \right\rangle.
\end{equation*}
\end{thm}

Table~\ref{tab:kjreg} offers the initial counting sequence for some
small values of $k$ and $j$. Each one corresponds to a known
differential equation satisfied by it generating function. In
Table~\ref{tab:graph} we compute the asymptotic number of some of
these graphs. 

We make one simple observation. The $\{k,k+2\}$-regular graphs are
isomorphic to the $k+2$-regular graphs with loops, by simply adding
loops to the vertices of degree $k$. This gives a family of identities
\begin{equation*}
  \left\langle E[e_2],\exp\left(t(h_i+h_{i+2})\right)\right\rangle
 =\left\langle E[h_2], \exp(th_{i+2})\right\rangle.
\end{equation*}

\begin{table}
\center
\begin{tabular}{cl}
$i,j$ & Initial terms in the counting sequence\\
1,2 &
1,0,1,4,18,112,820,6912,66178,708256,8372754,108306280,1521077404
\\
1,3& 
1, 0, 1, 0, 8, 0, 730, 0, 188790, 0, 102737670, 0, 102172297920,0
\\
1,4&
1,0,1,0,3,6,30,1011,38920,1920348,116400186,8580463110,757574641296
\\
2,3 &
1,0,0,1,10,112,1760,35150,848932,24243520,805036704,30649435140
\\
2,4&
1,0,0,1,3,38,730,20670,781578,37885204,2289786624,168879532980
\\
3,4&
1,0,0,0,1,26,820,35150,1944530,133948836,11234051976,1127512146540
\\
\end{tabular}\\

\caption{Counting sequences for $\{i,j\}$-regular graphs for small
  values of~$i$ and~$j$.}
\label{tab:kjreg}
\hrule
\end{table}

\subsection{Set covers and uniform hypergraphs} We now illustrate the
method on another family of objects, which results in set covers and
uniform hypergraphs.
An $n$-set is a set of cardinality~$n$.

\begin{defn}[$k$-cover of a set]
  A collection of $r$-sets $\mathcal{B}=\{B_1, \ldots, B_r\}$ is an {\em
  $r$-cover} of~$S$ if~$\bigcup_{i=1}^r B_i =S$. If~$S=[n]=\{1,2,\dots, n\}$, then it is an $r$-cover of $n$.  A cover is {\em
  restrictive\/} if all of the $B_i$ are distinct.  An $r$-cover is
  {\em $k$-regular} if any given element occurs in exactly $k$-subsets.
\end{defn}

A combinatorial argument shows that the number of distinct covers for
a set of $n$ elements is
\begin{equation*}
  \frac12\sum_{k=0}^n(-1)^k\binom{n}{k}2^{2^n-k},
\end{equation*}
which is clearly not P-recursive (equivalently, its generating series
is not D-finite.)

Devitt and Jackson~\cite{DeJa82} give a generating function for the
number of $k$-regular $r$-covers of $[n]$, a notion introduced by
Comtet~\cite{Comtet68}. Further, they prove that the number of
arithmetic operations required to actually calculate the number of
$k$-covers of an $n$ set by their method is bounded by $cn^k
\log n$. Results for fixed $k$, specifically $k=2, 3$ were treated by Comtet~\cite{Comtet68} and Bender~\cite{Bender74a} respectively.

We can derive enumeration formulas. For example, a $k$-regular graph
on~$n$ vertices is a restrictive $k$-cover of $[n]$ into 2-sets. In
general, calculating the generating function for restrictive
$k$-covers of $[n]$ into $j$-sets can be expressed as
\[
\rsc{\Gamma_{\sE\circ\sE_j}(p_1, p_2, \ldots)}{\sum_n h_{k}^nt^n}=
\rsc{\EE[e_j]}{\sum_n h_{k}^nt^n}.
\]
To determine $k$-covers with mixed-cardinality sets, say both $i$ and
$j$, we calculate
\[\rsc{\Gamma_{\sE\circ(\sE_i+\sE_j)}(p_1, p_2,
  \ldots)}{\sum_n h_{k}^nt^n}= \rsc{\EE[e_i+e_j]}{\sum_n h_{k}^nt^n}.\]

This yields the following simple consequence of
Theorem~\ref{thm:Dfin_species}. 
\begin{cor}
  Let $S$ be a finite set of integers.  For fixed $n$, and fixed $k$, the exponential generating
  function for $k$-regular $S$-covers of sets is D-finite, and is given by the scalar product
\begin{equation*}
\left\langle E[\sum_{s\in S} e_s], \exp(h_kt)\right\rangle.
\end{equation*}
\end{cor}

\begin{example}
We can express the problem of counting distinct
restrictive 2-covers of a set of cardinality $n$ by sets of
cardinality less than 5 as a scalar product. Denote the exponential generating
function of such set covers, by $S(t)$. We have,
\begin{equation*}
S(t)=\rsc{E[e_1+e_2+e_3+e_4]}{\exp(t h_2)}.
\end{equation*}

This problem is perfectly suited to either of our algorithms. We can
determine this differential equation, and the initial terms of the
counting sequence:
\[ 1, 0, 1, 8, 80, 1037, 17200, 350682, 8544641, 243758420, 8010360039.\]
\end{example}

It is worthwhile to remark that for a fixed~$j$, the set coverings by $j$-sets are equivalent to loopless $j$-uniform hypergraphs without multiplicities. These are encoded by $E[e_j]$. If we wish to encode hypergraphs with loops, we replace $e_j$ by $h_j$, and if we wish to encode hypergraphs with multiplicities we replace $E$ by $H$. 


\section{Asymptotic analysis} \label{sec:asympt}
Now that we have established how to determine the differential
equations satisfied by regular families of combinatorial objects, we
process these differential equations to obtain asymptotic enumeration
results.

Asymptotic enumeration of regular graphs is a topic that has received a great deal of attention. Indeed, as Gropp~\cite{Gropp92}
points out, the basic problem of regular graph enumeration was
considered before graphs were even ``invented'', over 120 years ago.
We first see some explicit results for graphs with fixed degree
sequences in the work of Read~\cite{Read58,Read60}, however, these are
rumoured to be ``difficult to penetrate''. Nonetheless, one can
determine an asymptotic expression for the number of 3-regular
graphs. Bender and Caufield~\cite{BeCa78} produce the first general
asymptotic formula for  the number of $k$-regular graphs on $n$
vertices, and Bollob\'as produces a similar result by a more
probabilistic approach that generalizes with ease to treat
hypergraphs. Next,  work by McKay~\cite{McKay90} and McKay and
Wormald~\cite{McWo91} consider the problem of $k$ which is not fixed,
but rather a function of $n$, and they achieve a formula which they
believe to be true in general, valid as $n\to\infty \text{ uniformly for } 1\leq k=o(n\sp {1/2})$
\begin{equation}
g_{k}(n)\sim \frac {(nk)!}{(nk/2)!2\sp {nk/2}(k!)\sp n}\exp\bigg(-\frac {k\sp
  2-1}4-\frac {k\sp 3}{12n}+O(k\sp 2/n)\bigg).
\end{equation}

This resembles Bollob\'as' asymptotic formula~\cite{Bollobas80} for labelled $k$-regular $r$-uniform hypergraphs, on $n$ vertices when $\frac{nk}{r}$ is an integer
\begin{equation*}
g_{k}^{(r)}(n)\sim\frac{(nk)!}{(nk/r)!(r!)^{(nk/r)}(k!)^n}\exp\left(-(r-1)(k-1)/2\right).
\end{equation*}
He also gives a formula for hypergraphs in which hyperedges only have single vertex intersections, which gives the constant in the McKay and Wormald formula for $r=2$. 

The asymptotic enumeration problem of regular graphs has been treated
with a variety of methods, such as the multidimensional Cauchy
integral technique mentioned earlier~\cite{McKay90}, a ``switching''
technique based on inclusion exclusion~\cite{McWa03,McWo91}, and some direct
combinatorial arguments on the equivalent problem of symmetric (0,1)
matrices with fixed row sum~\cite{BeCa78}.

Bollob\'as~\cite{Bollobas82} remarks that the number of
$k$-regular unlabelled graphs grows asymptotically like~$l_n/n!$ as~$n$
tends to infinity and where $l_n$ is the number of labelled regular
graphs. Intuitively, this is due to the fact that, for most large
graphs with no isolated vertices, and at most one vertex of maximal
degree, the automorphism group consists of only the identity automorphism.

The enumeration of other configurations is relatively untreated.
Gessel remarked~\cite{Gessel90} that the exponential generating functions of $k$-regular $r$-uniform hypergraphs (with and without loops, with and without multiplicities) are D-finite and the differential equations they satisfy are obtainable via the scalar
product. Domoco{\c{s}}~\cite{Domocos96} determines a scalar product form for
the generating series of minimal coverings that are multipartite
hypergraphs.

Here we continue to treat a variety of configurations. The results are
tabulated in Table~\ref{tab:graph} and Table~\ref{tab:other} and allow for
a comparison across objects rather than regularity parameter. All of the results were automatically generated. 

\subsection{Technique}
Our method is a classical singularity analysis of formal solutions of the linear
differential equations. It is precisely the same method we used in our
analysis of $k$-uniform Young tableaux~\cite{ChMiSa05}, and
thus we do not repeat the details here. Instead, after a short
description of the major steps, we present the fruits of our
analyses. A Maple worksheet of the computations is available at\\ {\tt
  http://www.math.sfu.ca/\verb+~+mmishna}.  

In the simplest cases, essentially the cases we could analyze directly with combinatorial arguments, we can solve the differential equation and do an asymptotic analysis on the solution.  In the more complex cases, we first convert our differential equation to the recurrence satisfied by the coefficients. Our series are D-finite, and thus such a sequence is bounded by a rational power of $n!$, and thus, we scale our sequence until it is convergent, and this allows us a more precise analysis. We convert this recurrence back to a differential equation, and determine the roots of the polynomial which is the coefficient of the leading term. From this we can calculate the dominant singularity, and determine a power series solution to the differential equation around this point. From this, we analyze the solution to determine an asymptotic expression  for the coefficients. This can be done automatically using tools from Maple, specifically DEtools and gfun. Finally, by generating sufficiently many terms in the sequence, we compare with the formula to determine a value for the constant.  

In the tables that follow, the Sloane sequence number refer to the
counting sequence as indexed in the Sloane On-line Encyclopedia of
Integer Sequences~\cite{Sloane}. 

Furthermore, we mention only the dominant term of the asymptotic
expression, but we could get the subsequent terms, save for the
appropriate constants. This is a consequence of the principal weakness
of our method-- we cannot generate exact expressions for the
constants.

We remark that, as presented in Table~\ref{tab:graph}, graphs of different types have the same asymptotic development, but differ only in the constants. 
\begin{table}
\center
\begin{tabular}{|p{3cm}||l|l|l|l|}\hline
$k$ & & 1 & 2& 3\\\hline\hline
Formula  &  \begin{minipage}{1cm}{\mbox{}\\ \vspace{1.5cm}}
              \end{minipage}
           &   $\displaystyle\left(\frac12\right)^{\frac n2}\frac{n!}{(n/2)!}$  
           & $\displaystyle\frac{n!}{\sqrt{n}}$
           & $\displaystyle\left(\frac32\right)^{\frac n2}\frac{n! (n/2)!}{n}$\\
\hline\hline
Simple graphs& Sloane \#  & A001147 & A001205 & A002829\\ 
$E[e_2]$& Constant & 1       & $e^{-\frac34}/\sqrt{\pi}$ & .043\\
\hline
Graphs with loops  & Sloane \#   & A001147 & A108246 & A110039\\
$E[h_2]$& Constant & 1       &  $e^{-\frac14}/\sqrt{\pi}$ & .318 \\ 
\hline
Multigraphs & Sloane \#   & A001147 & A002137 & A108243\\ 
$H[e_2]$& Constant & 1       & $e^{\frac14}/\sqrt{\pi}$ & .318\\
\hline
Multigraphs with loops & Sloane \#   & A001147 & A002135 & A005814\\ 
$H[h_2]$& Constant & 1       & $e^{\frac34}/\sqrt{\pi}$ & 2.35\\
\hline
\end{tabular}\\

\caption{
Asymptotic enumeration formulas for different classes of
  $k$-regular graphs. Formulas for 1- and 3- regular are valid only
  for even $n$.}

\label{tab:graph}
\hrule
\end{table}
\begin{table}
\center
\begin{tabular}{|l|c|l|l|l|}\hline
$k$ or $S$& Restrictions & Sloane \# & Formula & Constant \\\hline
\hline
\multicolumn{5}{|c|}{$S$-regular graphs}\\\hline
$\{1,2\}$ & & A00986 & $\displaystyle n^{-\frac 12}e^{\sqrt{2n}}n!$ &
$e^\frac{-3}2/2\sqrt \pi$\\
$\{2,3\}$ && A110040 & $\displaystyle n^{-\frac34} \left(\frac{\sqrt
  3}{2}\right)^n e^{\sqrt{3n}} n!^{3/2}$ & 0.007 \\
$\{1,3\}$ & $n=0\mod 2$ & A110039 & $\displaystyle n^{-1}\left(\frac32\right)^{n/2}n! (n/2)!$  & 0.43\\
$\{1,2,3\}$ & & A110041 &  $\displaystyle n^{-\frac
  34}\left(\frac{\sqrt 3}2\right)^n e^{\sqrt{3n}}n!^{3/2} $& 0.05\\
\hline
\multicolumn{5}{|c|}{$k$-regular 3-uniform hyper-graphs: $E[e_3]$}\\\hline
1 & $n=0\mod 3$ &A025035& $\displaystyle \left(\frac1{3!}\right)^{n/3}\frac{n!}{(n/3)!}$ & 1 \\
2 & $n=0\mod 3$ &A110100& $\displaystyle n^{-1}\left(\frac32\right)^{n/3}n!(n/3)! $ & 0.175\\
3 &   &A110101& $\displaystyle n^{-1}\left(\frac34\right)^nn!^2$ & 0.037\\
\hline
\multicolumn{5}{|c|}{$k$-regular 4-uniform hyper-graphs: $E[e_4]$}\\\hline
1 &$n=0 \mod 4$& A110102&$\displaystyle \left(\frac 1{4!} \right)^{n/4} \frac{n!}{(n/4)!}$ & 1\\
2 &$n=0 \mod 2$&A110103&$\displaystyle n^{-1}\left(\frac 23\right)^{n/2} n! (n/2)!$&0.100 \\
\hline
\end{tabular}

\smallskip

\caption{Asymptotic enumeration of different classes of regular objects}
\label{tab:other}

\hrule
\end{table}
There are a few observations we can make based on this data. We see
that allowing repetitions influences only the constant of the asymptotic
expansion, and often only slightly. We see this again in the cycle
covers in Table~\ref{tab:other}. The formulas look different than the
graph formulas we presented earlier, however if we expand the
factorials with Stirling's formula for $n!$, we quickly see that they
are the same. 

Although we are able to compute the differential equations for the
generating functions of the classes of 4-regular (multi-) graphs (with
loops), their asymptotic analysis is more complicated to do in an
automated fashion, because a saddle point analysis arises. This is an obvious starting
point for future work. Most of the tools are already implemented,
it is mostly a question of understanding them, and determining how to
best automate them.

We could make further observations by considering directed versions of
any of these structures. For directed graphs, we need only
consider $\sE\circ\sL_2$, where $\sL_k$ is the species of lists of
length $k$, and we could generalize hypergraphs in a number of ways;
by putting an order, or even an orientation on each ``edge'' using the
species of cycles or lists, as in the directed graph case. 

\section{Comments, conclusions and perspectives}
\subsection{Asymptotic expansions of different families of functions}
Coefficients of taylor expansions of algebraic functions have a known kind of expansion, that can in fact we used to establish the transcendence of some series~\cite{Flajolet87}. We can also describe asymptotic discrepancy criteria for coefficients of D-finite functions. 

As we remarked earlier, coefficients of D-finite series are also restricted in their asymptotic growth. A more complete version of this criteria is following theorem presented in Wimp and Zeilberger~\cite{WiZe85}.
\begin{thm}[Wimp and Zeilberger] Suppose that $f(t)=\sum_{n\geq 0} f_n t^n$ is a
  D-finite series in in $\C[[t]]$. Then, for sufficiently large $n$,
  the coefficients $f_n$ have an asymptotic expansion which is a sum of terms of the form
\[\lambda(n!)^{r/s}\exp(Q(n^{1/m}))\omega^nn^\alpha(\log n)^k,\] 
  where $r, s, m, k\in\N$, $Q$ is a polynomial and $\lambda,
  \omega, \alpha$, are complex numbers.
\end{thm}
We may ask ourselves, have our examples encompassed the full asymptotic potential of D-finite sequences?  We are very curious about the combinatorial structure of
families which do have such expansions. Are D-finite
species sufficient to consider?

For example, in our earlier study of $k$-uniform Young
Tableaux~\cite{ChMiSa05}, which are enumerated by the scalar product
$\left\langle H[e_1+e_2], \exp(h_kt)\right\rangle$ have a conjectured
form (verified for $k=1..4$) of
\begin{equation*}
y_n^{[k]}\sim
\frac1{\sqrt2}\left(\frac{e^{k-2}}{2\pi}\right)^{k/4}
n!^{k/2-1}\left(\frac{k^{k/2}}{k!}\right)^n\frac{\exp(\sqrt{kn})}{n^{k/4}},
\qquad n\rightarrow\infty.
\end{equation*}
This is an exact conjecture, more complete than the examples we have presented here,
although it results from the same kind of calculation, and presumably if we
completed the complex saddle point analyses required for the
$k=4$ cases, we might be able to guess such a form for our examples.    

It is also of interest to note that while in some cases
operations of summation and integration preserve D-finiteness, $G_r$,
the class of {\em all\/} regular graphs is not D-finite. The same is
true of all the classes we have presented here: Although for any $k$,
the subclass of $k$-regular objects is D-finite, the larger subclass
of regular objects is not.
This is interesting and can help us refine our notion of D-finiteness. 

\subsection{D-finite species?} \label{sec:conclusion}
We have thus far restrained ourselves from defining a notion of
D-finite species. Ideally, such a theory would contain two main
components: A ``D-finite  species'' should satisfy some sort of system combinatorial
differential equations with polynomial coefficients; and symmetric
series of such species should be D-finite. We would then expect to be
able to have theorems of the form:
\begin{enumerate}
\item\label{eq:cis1} If \sF and \sG are D-finite species, then so are
     $\sF+\sG$, and $\sF\cdot \sG$, $\sF'$;
\item\label{eq:cis2} If \sG is a polynomial species, (in particular,
if its cycle index series is a polynomial), then $\sF\circ\sG$ is a
D-finite species;
\item\label{eq:cis3} If \sF satisfies an ``algebraic equation'' of
  species, including for instance, equations of the form
  $\sF=\sX\mathsf{P}(\sX, \sF)$ for polynomial species $\mathsf{P}$,
then $\sF$ is a D-finite species.
\end{enumerate}
A candidate definition is given in~\cite{Mishna03}, however more work
remains to be done. 

Finally, much work has been done to characterize combinatorial classes
of objects with rational and algebraic generating series
(see~\cite{MBM06} for a recent summary), and
hopefully this work is a step towards such a characterization for
D-finite generating functions. We are encouraged by a recent thought
of  Flajolet, Gerhold and Salvy~\cite{FlGeSa05},
\begin{quote}
\em Almost every thing is non-holonomic unless it is holonomic by design.\\
\end{quote}
(Series are holonomic if and only if they are D-finite.) They follow
this with the remark that there are several surprising exceptions to
this rule, notably $k$-regular graphs. Hopefully we have demonstrated
that this isn't so surprising; That in fact, there are deep reasons
underlying the D-finiteness of objects with this sort of regularity,
and that furthermore, this D-finiteness can be exploited in an
automatic way.

\subsubsection*{Acknowledgements} The author wishes to thank Fran\c cois
Bergeron, Bruno Salvy, Fr\'ederic Chyzak, and indirectly Fr\'ederic
Jouhet, for fruitful, instructive discussions. Thank you also to Cedric
Chauve for useful comments on the text.

Sequences:  A001147, A001205, A002829,  A108246, A110039,  A002137,
A108243,  A002135, A005814,A00986 A110040, A110041, A025035 A110100,A110101, A110102, A110103.
\bibliographystyle{acm}

\end{document}